\documentclass[12pt,a4paper]{article}
\usepackage[T2A]{fontenc}
\usepackage[cp1251]{inputenc}
\usepackage{amsmath,amssymb,amsthm,amsfonts,amscd}

 \theoremstyle{plain}
\newtheorem{theorem}{Теорема}
\newtheorem*{main*}{Основная Теорема}
\newtheorem*{theorem*}{Теорема}
\newtheorem{lemma}[theorem]{Лемма}
\newtheorem{proposition}[theorem]{Предложение}
\newtheorem{corollary}[theorem]{Следствие}
\newtheorem{definition}[theorem]{Definition}

\newtheorem*{remark*}{Remark}
\newtheorem*{example*}{Example}
\newtheorem{example}[theorem]{Example}

\def\Z{{\Bbb Z}}
\def\R{{\Bbb R}}
\def\RP{{\Bbb R}\!{\rm P}}

\def\dim{{\operatorname{dim}}}

\def\int{\operatorname{\rtimes}}

\def\Z{{\Bbb Z}}
\def\R{{\Bbb R}}
\def\RP{{\Bbb R}\!{\rm P}}

\def\R{{\Bbb R}}

\begin{document}

\sloppy
\date{}

\title{A geometrical approach toward stable homotopy groups of spheres. A Desuspension Theorem (revisited)}

\author{Petr Akhmetiev}

\maketitle

\begin{abstract}
Desuspension Theorem for stable steam $\Pi_{2^l-2}$ in stable homotopy groups of spheres is formulated and is proved. This result is a minor revision of a theorem in \cite{A1}.
\end{abstract}

\section{Introduction}

A Desuspension Problem in stable homotopy groups of spheres is following.
Let $\Pi_n$ be the stable homotopy group of sphere of dimension $n$, 
it is well-known that $\Pi_n=\pi_{2n+2}(S^n)$. The following homomorphism, which is called  the $k+2$-desuspension  is well-defined: 
$$\Sigma^{k+2}: \pi_{2n-k}(S^{n-k-2}) \to \pi_{2n+2}(S^n), \quad k\ge -1.$$
What is the maximal number
 $d(n)$, for which the homorphism $\Sigma^{d+2}$ is an epimorphism? 
 The Freudenthal Desuspension Theorem clams
  $d(n) \ge -1$ for an arbitrary $n \ge 1$.  By the Adams Theorem on Hopf Invariant for
 $n=1,3,7$ one gets: $d(n)=1$, and for $n \ne 0,1,3,7$ one gets: $d(n) \ge 2$.
 Oppositely, for the Mahowald element
 $\eta_l \in \Pi_{2^l}$, $l \ge 3$ one gets: $d(2^l)=4$. Our main result is the following: 

\begin{theorem}{Desuspension Theorem}\label{comp}

	For an arbitrary positive integer $d$ there exists a positive
	integer  $l_0=l_0(d)$ such that an arbitrary element in the stable homotopy group of spheres
	$\Pi_{2^l-2}$,  $l \ge l_0$, admits a desuspension
	of the order $d+2$.
	\end{theorem}

\section{Compression Theorem}

Let us consider an approach toward stable homotopy groups of spheres, based on the Pontryagin approach \cite{P} and the Smale-Hirsh Theorem on immersions \cite{H}. Results using the considered approach was previously developed, in particular, in \cite{C1}, \cite{C2}. By this approach the stable homotopy group of spheres
 $\Pi_n$ is naturally isomorphic to the cobordism group $Imm^{fr}(n,1)$ of skew-framed immersions in the codimension $1$. An element of this cobordism group is represented by an immersion   $f:N^n \looparrowright \R^{n+1}$, where the manifold  $N^n$ is oriented. The group structure is well-defined by disjoint union of immersions, a cobordism class of immersion is defined as the regular cobordism class.

Let us define the cobordism group
 $Imm^{sf}(n-1,1)$ as equivalent classes of immersions
$g: M^{n-1} \looparrowright \R^n$, where, instead of $Imm^{fr}(n-1,1)$,  the manifold $M^{n-1}$ is non-oriented. 

\begin{definition}\label{opr1}
Let us call an element	
 $[g] \in Imm^{sf}(n-1,1)$ admits a compression of an order  $d$, $d \le n-1$, if there exists an immersion   $g_1: M^{n-1}_1 \looparrowright \R^n$, which is cobordant to the immersion $g$, 
 for which the orienting cohomology class  (this class is denoted by  $\kappa_1 \in H^1(M^{n-1}_1))$) is represented by a mapping
 $\kappa_1: M_1^{n-1} \to \RP^{\infty}$ with the image into the standard skeleton
 $\RP^{n-1-d} \subset \RP^{\infty}$.
\end{definition}

Let
$$ \lambda: Imm^{sf}(n-1,1) \to Imm^{fr}(n,1)$$
be the Khan-Priddy transfer, which is formulated for cobordism group and is called the Koschorke homomorphism. It is well-known that the Koschorke homomorphism is a slitted epimorphism onto $2$-component. 

\subsection*{Compression Problem}

What is the maximal integer $d(n)$, for which an arbitrary element $[g] \in Imm^{sf}(n-1,1)$
admits a compression of the order $d$?

\subsection*{Desuspension Problem}

What is the maximal integer
 $d(n)$, for which an arbitrary element $[f] \in Imm^{fr}(n,1)$ admits a desuspension of the order  $d+2$?
\[  \]

The Compression Problem is equivalent to the Desuspension Problem in  metastable dimensions $2d + 1 \le n$. Therefore the Desuspension Theorem  \ref{comp} is a corollary of the following result.

\begin{theorem}{Compression Theorem}\label{comp1}

For an arbitrary positive integer $d$ there exists a positive
integer  $l=l(d)$ such that an arbitrary element in the cobordism
group $Imm^{sf}(2^{l'}-3,1)$, for $l' \ge l$, admits a compression
of the order $d-1$.
\end{theorem}


\[  \]

\section{The cobordism group $Imm^{sf}(n-k,k)$, $k \ge 1$ and its generalization}

Let $M^{n-k}$ be a closed manifold of dimension $(n-k)$, $n \ge k \ge 1$, $\varphi:
M^{n-k} \looparrowright \R^n$ be an immersion of this manifold
into $\R^n$ in the codimension $k$, $\Xi_M$ is a skew-framing of
the immersion $\varphi$ with the characteristic class $\kappa_M$
of this skew-framing. Additionally, let us assume that the
manifold $M^{n-k}$ is equipped with the family of 1-dimensional
cohomology classes modulo 2:
\begin{eqnarray}\label{A_M}
A_j=\{\kappa_i\}, \quad \kappa_i \in H^1(M;\Z/2), \quad 0 \le i
\le j, \quad \kappa_0 = \kappa_M.
\end{eqnarray}
This collection of cohomology classes is represented by the
collection of classifying maps:
\begin{eqnarray}\label{A'_M}
A_{j}=\{\kappa_i: M^{n-k} \to \RP^{\infty} \}, \quad i=0,1,
\dots, j, \quad \kappa_0 = \kappa_M.
\end{eqnarray}

\begin{definition}\label{cob}

The cobordism group $Imm^{sf;A}(n-k;k)$ is
represented by triples $(\varphi,\Xi_M,A_j)$, where:

 ---$\varphi: M^{n-k} \looparrowright \R^n$--is an immersion of a closed $(n-k)$-dimensional manifold into Euclidean space,

--- $\Xi_M$ is skew-framing of the immersion $\varphi$,

---$A_j$ is a collection  of cohomology classes, described in $(\ref{A'_M})$.

The cobordism relation of triples  is the standard.
\end{definition}

\subsubsection*{Remark}
In the case $j=0$ the cobordism group $Imm^{sf;A_j}(n-k;k)$
coincides with the cobordism group $Imm^{sf}(n-k;k)$, which is called the cobordism group of
skew-framed immersions.
\[  \]

A natural homomorphism
\begin{eqnarray}\label{102}
J^{k}_{sf;A_j}: Imm^{sf;A_j} (n-1,1) \to Imm^{sf;A_j}(n-k,k)
\end{eqnarray}
is defined as follows. Let us assume that the triple $(\varphi_0,
\Xi_{M_0}, A_j(M_0))$ represents an element in the cobordism group
 $Imm^{sf;A_j}(n-1,1)$. Let us consider the following triple $(\varphi, \Xi_M, A_j(M))$,
where the immersion $\varphi: M^{n-k} \looparrowright \R^n$ define
as follows. The manifold $M^{n-k}$ is a submanifold in the
manifold $M_0^{n-1}$, the fundamental class of this submanifold
represents the homological Euler class of the bundle
$(k-1)\kappa_{M_0}$, the immersion $\varphi$ is defined as the
restriction of the immersion $\varphi_0$ on $M^{n-k}$. The
skew-framing $\Xi_M$ of the immersion $\varphi$ is defined by the
standard construction like in the case $j=0$, the collection
$A_j(M)$ of cohomology classes is the restriction of the
collection $A_j(M_0)$ on the submanifold $M^{n-k} \subset
M_0^{n-1}$.
\[  \]

Let us generalize Definition $\ref{opr1}$ for the cobordism group
$Imm^{sf;A_j}(n-k,k)$.

\begin{definition}\label{comp-5}
Let $[(\varphi,\Xi_M,A_j)] \in  Imm^{sf;A_j}(n-k,k)$. We shall say
that the element  $[(\varphi,\Xi_M,A_j)]$ admits a compression of
the order  $d$, if in its cobordism class there exists a triple
 $(\varphi',\Xi'_M,A'_j)$, such that the pair  $(M^{n-k}, \kappa_0)$
admits a compression of the order $d$ is the sense of the
Definition $\ref{opr1}$.
\end{definition}

Let us define the transfer homomorphism
\begin{eqnarray}\label{transfer}
r_j^!:Imm^{sf;A_j}(n-k,k) \to Imm^{sf;\hat A_{j-1}}(n-k,k)
\end{eqnarray}
 with respect to the cohomology
class $\kappa_j$.

Let $x \in Imm^{sf;A_j}(n-k,k)$ be an element represented by a
triple $(\varphi, \Xi_M, A_j)$. Let us define the 2-sheeted cover
$$p_j:\hat M^{n-k} \to M^{n-k}$$
as the regular cover with the characteristic class
 $w_1(p_j)=\kappa_j \otimes \kappa_M \in H^1(M^{n-k};\Z/2)$.
(We will denote below by $p_j$ the linear bundle over $M^{n-k}$
with the characteristic class $w_1(p_j)$ and also the
characteristic class $w_1(p_j)$ itself.)

Let us define a skew-framing
 $\Xi_{\hat M}$. Let us consider the immersion
$\varphi: M^{n-k} \looparrowright \R^n$. Let us denote the
immersion $\varphi \circ p_j$ by $\hat \varphi$. Let us denote the
normal bundle of the immersion $\varphi$ by $\nu_{M}$. Let us
denote the normal bundle of the immersion $\varphi \circ p_j$ by
$\nu_{\hat M}$. Let us define the skew-framing $\Xi_{\hat M}$ of
the immersion $\varphi \circ p_j$ by the formula $\Xi_{\hat M} =
p_j^{\ast}(\Xi_M)$.

Let us define the collection of the cohomology classes $\hat
A_{j-1}$ by the following formula:
$$\hat A_{j-1} = \{ \hat \kappa_j = p_j \circ \kappa_j, \quad j=0, \dots, j-1 \}. $$

 We will define
$$r_j^!(x)=[(\hat \varphi,\Xi_{\hat M},\hat A_{j-1})].$$

\begin{example}\label{Eccles}
In the case  $j=0$ the transfer homomorphism $(\ref{transfer})$ is
given by the following formula:
$$r_{0}: Imm^{sf}(n-k,k) \to
Imm^{fr}(n-k,k)$$ The properties of this homomorphism is investigated in [A-E].
\end{example}

Let us consider the composition of the
$(j_2-j_1)$ transfer homomorphisms:
\begin{eqnarray}\label{103}
r_{j_1+1, \dots, j_2}^!: Imm^{sf;A_{j_2}}(n-k,k) \to Imm^{sf;\hat
A_{j_1}}(n-k,k), \quad j_2 \ge j_1.
\end{eqnarray}
In the case $j_1=0$, the homomorphism $(\ref{103})$, where $\psi=j_2$, is
denoted by
\begin{eqnarray}\label{103tot}
r_{tot}^!: Imm^{sf;A_{\psi}}(n-k,k) \to Imm^{sf}(n-k,k).
\end{eqnarray}

Let us describe the homomorphism $(\ref{103})$ explicitly. Let us
assume that an element $x \in Imm^{sf;A_{j_2}}(n-k,k)$ is given by
a triple $(\varphi, \Xi_M, A_{j_2})$. Let us consider the the
following subcollection $\{\kappa_{j_1+1}, \dots, \kappa_{j_2} \}$
of the last $(j_2-j_1)$ cohomology classes of the collection
$A_{j_2}$. Let us define the $2^{j_2-j_1}$--sheeted cover $p: \hat
M^{n-k} \to M^{n-k}$, given by the following collection of the
cohomology classes:
$$\{\kappa_{j_1+1} \otimes \kappa_0, \dots, \kappa_{j_2} \otimes \kappa_0 \}.$$
The immersion $\hat \varphi: \hat M^{n-k} \looparrowright \R^n$ is
given by the following composition:
$$ \hat \varphi = \varphi \circ p. $$
The normal bundle of the immersion  $\hat \varphi$ is equipped
with the skew-framing $p^{\ast}(\Xi_M)=\Xi_{\hat M}$. The
collection of cohomology classes
$$\hat A_{j_1} = \{ \hat \kappa_0
= \kappa_0 \circ p, \dots \hat \kappa_{j_1}=\kappa_{j_1} \circ p
\}$$
is well-defined by the classifying map.

The cobordism class of the triple $(\hat \varphi, \Xi_{\hat M},
\hat A_{j_1})$ determines the element $r_{j_1+1, \dots,
j_2}^!(x)$.

\begin{proposition}\label{totaltrans}
For arbitrary positive integers $l$ and $k$, 
$2^l-2=n>k>0$, there exists a positive integer $\psi = \psi(k,l)$,
such that the total transfer homomorphism $(\ref{103tot})$:
\begin{eqnarray}\label{104}
r_{tot}^!: Imm^{sf;A_{\psi}}(n-k,k) \to Imm^{sf}(n-k,k),
\end{eqnarray}
where $A_{\psi}$ is an index for  collections of $\psi$ classes, is trivial.
\end{proposition}


Let us start a proof of Proposition \ref{totaltrans} with the following lemma.

\begin{lemma}\label{P-Th}

The cobordism group $Imm^{sf;A_j}(n-k,k)$ is a finite $2$-group if
$n-k >0$, or, $n-k=0$ and $k$ is odd.
\end{lemma}

\subsubsection*{Proof of Lemma $\ref{P-Th}$}

The cobordism group $Imm^{sf}(n-k,k)$ by the Pontrjagin-Thom
construction is the stable homotopy group $\Pi_{n}(P_{k-1})$,
where $P_{k-1} = \RP^{\infty}/\RP^{k-1}$, see f.ex. [A-E]. By the
standard arguments for Thom spaces, the cobordism group $Imm^{sf;A_j}(n-k,k)$ is a
stable homotopy group $\Pi_{n}(P_{k-1} \times \prod_{i=1}^{j}
\RP^{\infty}_i/(pt \times \prod_{i=1}^{j}
\RP^{\infty}_i))$, $pt \in P_{k-1}$ is a marked point. The space $(P_{k-1} \times \prod_{i=1}^{j} \RP^{\infty}_i/(pt \times \prod_{i=1}^{j}
\RP^{\infty}_i))$ is an $(k-1)$-connected  2-homotopy type space (unless $n = k \equiv 0 \pmod{2}$). Lemma
$\ref{P-Th}$ is proved. \qed

Let us define the following sequence of positive integers
 $s_{n-k}, \dots,
s_1$, where the indexes decrease from $(n-k)$ to 1:
$$ s_{n-k} =
ord(\Pi_{n-k}), \quad s_{n-k-1}= ord(\Pi_{n-k-1}), \quad \dots,
\quad s_1= ord(\Pi_1).$$
 In this formula we denote by
$ord(\Pi_i)$ the logarithm of the maximal order of an element in
the 2-component of the $i$-th stable homotopy group of spheres.
Then let us define
\begin{eqnarray}\label{psi1}
 \psi(i)=1+\sum_{j=i}^{n-k} s_j
 \end{eqnarray}
 and let us define the required integer
 by the following formula:
 \begin{eqnarray}\label{psi}
\psi = \psi(1)+ \sigma,
\end{eqnarray}
where $\sigma = ord(Imm^{sf;A_{\psi(1)}}(n-k,k))$ the logarithm of the order of the cobordism group, where the number $\psi(1)$ is defined by (\ref{psi1}) using Lemma
$\ref{P-Th}$. Therefore, we have $\psi \ge \psi(1) \ge \psi(2) \ge
\dots \ge \psi(n-k) \ge 1$.

Let  $x \in Imm^{sf;A_{\psi}}(n-k,k)$ be an  element represented
by a triple $(\varphi, \Xi_M, A_{\psi})$. Let us consider the
element
\begin{eqnarray}\label{element}
r^!_{\psi - \psi_1}(x) = (\hat{\varphi},\hat{\Xi},\hat{A}), \quad \hat{A}_{\psi_1}=\{\hat \kappa_{\psi_1}, \dots, \hat \kappa_{\psi_{n-k}},\hat{\kappa}_M\},
\end{eqnarray}
and the element
\begin{eqnarray}\label{element1}
r^{forg}_{\psi - \psi_1}(x) = (\varphi, \Xi, A^{forg}_{\psi_1}), \quad A^{forg}_{\psi_1}=\{\kappa_{\psi_1}, \dots, \kappa_{\psi_{n-k}},\kappa_M \}.
\end{eqnarray}
Replies the element (\ref{element}) by the following triple:
\begin{eqnarray}\label{element2}
r^!_{\psi - \psi_1}(x) \mapsto r^!_{\psi - \psi_1}(x) + \sigma r^{forg}_{\psi - \psi_1}(x),
\end{eqnarray}
where $\sigma$ is given by  (\ref{psi}). 
Let us prove that the element $(\ref{element2})$
admits a compression of the order 0, see Definition
$\ref{comp-5}$.

Let $(\varphi, \Xi_M, A_{\psi})$, $\varphi: M^{n-k} \looparrowright \R^n$, be the given triple.
Denote the product  $\RP^{\infty}(0) \times
\prod_{j=\psi_1}^{j=\psi} \RP^{\infty}(j)$ by   $X(\psi_1)$, and let us
consider the map
\begin{eqnarray}\label{lambdaM}
\lambda'_M: M^{n-k} \to X(\psi_1),
\end{eqnarray}
defined as the
direct product of  the classifying maps of the subcollection
$A^{forg}_{\psi_1} \subset A_{\psi}$ of cohomology classes, given by (\ref{element1}). Let us
denote the space $X(\psi_1)$ by $X$ for short.

Let us consider a natural filtration
\begin{eqnarray}\label{105}
\dots \subset X^{(n-k+1)} \subset X^{(n-k)} \subset \dots \subset
X^{(1)} \subset X.
\end{eqnarray}
This filtration is the direct product of the standard
coordinate filtrations. Each stratum  $X^{(i)} \setminus
X^{(i+1)}$, $i=1, \dots n-k$ is an union of open cells of the
codimension $i$. Each cell is determined by the corresponded
multi-index $\mu=(m_0,m_{\psi_1} \dots, m_{\psi})$, $m_0+ m_{\psi_1}+ \dots +
m_{\psi}=i$ is a deep of a strata, each lower index in a multi-index
shows a corresponded
coordinate projective space that contains the given cell.

Let us assume that the map $\lambda$ is in a general position with
respect to the filtration  $(\ref{105})$. Let us denote by $L^0
\subset M^{n-k}$ a $0$--dimensional submanifold in $M^{n-k}$,
defined as the inverse image of the total stratum  $X^{(n-k)}$ of the deep $n-k$ of the
filtration.

Let us consider the triple $(\hat \varphi, \hat{\Xi}_{\hat M}, \hat
A_{\psi(1)})$, given by (\ref{element}). Let us prove that this triple is cobordant to
a triple $(\hat \varphi', \hat\Xi'_{\hat M'}, \hat A'_{\psi(1)})$
such that the mapping
 $\hat \lambda': \hat M'^{n-k} \to X$,
 constructed by means of the collection of the cohomology classes
 $A_{\hat M'}$, satisfies the following property:
 \begin{eqnarray}\label{106}
  (\hat \lambda')^{-1}(X^{(n-k)}) = \emptyset.
\end{eqnarray}

Let us consider arbitrary
collection of orientations of cells in $X^{(n-k)} \setminus
X^{(n-k+1)}$. Let us take a local orientation of the bundle $\kappa_M$, this local orientation is defined by the corresponding local orientation of the universal bundle $\gamma$, where $\kappa_M$ is a pull-back of $\gamma$ over the component $\RP^{\infty}(0)$ by $\lambda$. 
In particular, one may speak on  orientation of $\kappa_M$ over each cell of $X$, which is defined when we fix one of the two sheets of the covering of $\kappa_M$. This common convention concerning signs  is considered
for $\hat{\lambda}$ and $\lambda$ simultaneously.

Let us  denote by $lk(\mu)$ the integer coefficient
of intersection of the image $\lambda(M^{n-k})$ with the
cell of a multi-index  $\mu$. Let us denote analogously
by $lk(\hat \mu)$ the integer coefficient of intersection of
the image $\hat \lambda(\hat M^{n-k})$ with the  cell of a prescribed
multi-index $\mu$. Obviously, the collection of the integers
$\{lk(\hat \mu)\}$ is obtain from the collection of the
corresponded  integers $\{lk(\mu)\}$ by the multiplication on
$2^{\psi-\psi(1)}$.

Let $A^{forg}_{\psi(1)}$ be the subcollection of $A_{\psi}$, consists of
the last $\psi(1)$ cohomology classes. By the construction, the
exponent of the group $Imm^{sf;A_{\psi(1)}}(n-k,k)$  equals to
$2^{\sigma}$. Therefore the disjoint union of the $2^{\sigma}$ copies of the triple $(\varphi, -\Xi_{M}, A^{forg}_{\psi(1)})$
(where by the opposition of the skew-framing we means the same triple with the inverse local orientation of skew-framing), 
determines the trivial element in the group
$Imm^{sf;A_{\psi(1)}}(n-k,k)$.

Let us consider the triple $(2^{\sigma})(\varphi, -\Xi_{M},
A_{\psi-\psi(1)})$. This triple is define as the disjoin union of
$2^{\sigma}$ copies of the triple  $(\varphi,
-\Xi_{M}, \hat{A}_{\psi(1)})$ (note, that the skew-framing $\Xi_M$ is changed
into the skew-framing with the opposite global orientation in the case $k$ is even). 
The collection of the coefficients for
the triple $(2^{\sigma})(\varphi, -\Xi_{M}, \hat{A}_{\psi-\psi(1)})$ is denoted by $\{lk(2^{\sigma}(-\mu))\}$. Obviously,
$\{lk(2^{\sigma}(-\mu))\} = - 2^{\sigma}\{lk(\mu)\}$.

Let us consider the triple  $(\hat \varphi, \hat \Xi_{\hat M}, \hat
A_{\psi(1)})$, and  the disjoint union of the triple
 with the triple
$(2^{\sigma})(\varphi, -\Xi_{M}, A^{forg}_{\psi(1)})$. The denote this new triple by
$(\hat \varphi_1, \hat{\Xi}_{1;\hat M_1}, \hat A_{1;\psi(1)})$. This new triple and the triple $(\hat \varphi, \hat \Xi_{\hat M}, \hat
A_{\psi(1)})$ represent a common
element in the cobordism group $Imm^{sf; A_{\psi(1)}}(k,n-k)$.
The mapping $\hat \lambda_1: \hat M_1^k \to X$, constructed by
means of the collection $\hat A_{1;\psi(1)}$ of the cohomology
classes is well-defined. The collection of the intersection
coefficients, defined for the mapping  $\hat \lambda_1 $ will be
denoted by $\hat lk(1;\mu)$. Obviously for an arbitrary multi-index
$\mu$ we have $\hat lk(1;\mu)=0$.

A normal surgery of the triple $(\hat \varphi'_1, \Xi_{1;\hat M_1},
\hat A_{1;\psi(1)})$  to a triple $(\hat \varphi_2, \Xi_{2;\hat M_2},
\hat A_{2;\psi(1)})$ by 1-handles is defined such that the the map
$\lambda_2$ is defined by means of the collection  $\hat
A_{2;\psi(1)}$ of cohomology classes, satisfies the condition
$$(\hat \lambda_2)^{-1}(X^{(n-k)}) =  \emptyset,$$
which is analogous to the condition $(\ref{106})$.
This gives the first step of the proof.

Let us describe  next steps of the proof. Let us denote the
triple $(\hat \varphi_2, \Xi_{2;M_2}, A_{2;\psi(1)})$ by $(\varphi,
\Xi_M, A_{\psi(1)})$ for short. Let us consider the subspace $X(\psi_2) \subset X(s)$, defined by the formula:
\begin{eqnarray}\label{X1}
X(\psi_2) = \RP^{\infty}(0) \times
\prod_{j=\psi_2}^{j=\psi} \RP^{\infty}(j),
\end{eqnarray}
\begin{eqnarray}\label{incl2}
	i_{\psi(2)}: X(\psi(2))
	\subset X(\psi(1)) \subset X(\psi). 
\end{eqnarray}
Denote $X(\psi_2)$ by $X$ again for short. The map $\lambda:  M^{n-k} \to X$,
constructed by means of the collection $A^{forg}_{\psi(2)}$ of cohomology
classes, satisfy a condition, which is analogous to $(\ref{106})$. Let us consider the
triple $(\hat \varphi, \Xi_{\hat M}, \hat A_{\psi(2)})$, given by
the following element
$$r^!_{\kappa_{\psi(1)}, \dots, \kappa_{\psi(2)-1}}(\varphi,
\Xi_M, A_{\psi(1)})$$ in the cobordism group $Imm^{sf; 
A_{\psi(2)}}(n-k,k)$.

Let us define the space $X(\psi(2))$, which is an analog of the space $X(\psi_1)$, as the Cartesian product of
infinite--dimensional projective spaces with indexes $(0,
\psi(2), \dots, \psi)$:
$$X(\psi(2)) = \RP^{\infty}_0 \times
\prod_{j=\psi(2)}^{j=\psi} \RP^{\infty}(j).$$ Denote $X(\psi(2))$ by $X$ for short.

 The space $X$ is equipped with the
following standard stratification (\ref{105}). 
The  inclusion $(\ref{incl2})$ is agree with the stratifications
$(\ref{105})$.

The collection $\hat{A}_{\psi(2)}$ of cohomology classes determines the
map  $\hat \lambda: \hat M^{n-k} \to X$. The condition
$(\ref{106})$ implies the following analogous condition for $\hat \lambda$:
\begin{eqnarray}\label{106bis}
\hat \lambda^{-1}(X^{(n-k)}) = \emptyset.
\end{eqnarray}
Let us denote by $\hat L^1 \subset \hat M^{n-k}$ a 1-dimensional
submanifold in $\hat M^{n-k}$ given by the following formula:
 $$\hat L^1= \hat \lambda^{-1}(X^{(n-k-1)}).$$
The restriction of the cohomology classes of the collection $\hat
\hat A_{\psi(2)}$ on the submanifold $\hat L^1$ is trivial. In
particular, the submanifold  $\hat L^1$ is framed.

The components of the manifold $\hat L^1$ are equipped with the
collection of the multi-indexes corresponded to the top cells of
the subspace   $X^{(n-k-1)}$. A fixed multi-index determines a
disjoint collection of  $2^{s_1}$ copies of 1-dimensional framed
manifold (probably, non-connected) and the copies are pairwise
diffeomorphic as a framed manifolds.

A framed 2-dimensional manifold $\tilde K^2$ with a framed
boundary $\partial (\tilde K^2)= (\tilde L^1)$ is well-defined.
This framed manifold determines the body of a handle for the
normal surgery of the triple $(\hat \varphi, \Xi_{\hat M}, \hat
A_{\psi(2)})$ to a triple $(\hat \varphi_1, \Xi_{1;\hat M_1}, \hat
A_{1;\psi(2)})$ such that the collection $\hat A_{1;\psi(2)}$ of
cohomology classes determines the map
 $$\hat \lambda_1: \hat
M_1^{n-k} \to X$$ satisfies the condition
\begin{eqnarray}\label{107bis}
	\hat \lambda_1^{-1}(X^{(n-k-1)}) = \emptyset.
\end{eqnarray}

The next steps of the proof are analogical to the step with the considered case with one-dimensional obstruction. The
parameter $i$ denoted the dimension of the obstruction is changed
from $2$ up to $n-k$. In each step  analogous conditions to the
conditions $(\ref{106})$, $(\ref{107bis})$ is considered. At the
last step of the proof we have a framed manifold
$(\hat{M}^{n-k},\hat{\Xi}_M)$, equipped with a collection $\hat A_{\psi(n-k)}$
of the trivial cohomology classes. The framed manifold
represented by $\psi(n-k)$ disjoint copies of the framed manifold
$(\hat{M}^{n-k},\hat{\Xi}_M)$ is a framed boundary (and therefore a
skew-framed boundary). Proposition $\ref{totaltrans}$ is
proved. \qed
\[  \]

Let us describe an algebraic obstruction for the compression of
a given order.

\begin{lemma}\label{codcompr}
An arbitrary element $x \in Imm^{sf;A_j}(n-k,k)$ admits a
compression of an order $i$, $i \le n-k$, if and only if the
element
$$J_{sf}^{k'}(x) \in  Imm^{sf;A_j}(n-k',k')$$
($i \le n-k'
\le n-k$) admits a compression of the same order $i$, where the homomorphism 
$J_{sf}^{k'}(x)$ is given by $(\ref{102})$.
\end{lemma}

\begin{corollary}\label{totalobstr}
For an arbitrary element $x \in Imm^{sf;A_j}(n-k,k)$ the total
obstruction for a compression of an order $q$  ($0 \le q \le n-k$)
is given by the element $J_{sf;A_j}^{q}(x) \in
Imm^{sf;A_j}(q,n-q)$.
\end{corollary}

To prove  Lemma $\ref{codcompr}$ and Corollary $\ref{totalobstr}$,
let us formulate an auxiliary lemma. Let as assume that a triple
$(\varphi, \Xi_M, A_j)$ represents an element
 $x \in Imm^{sf;A_j}(n-k,k)$. Let us additionally assume that
this element admits a compression of the order $(i-1)$. This means
that in the triple $(\varphi, \Xi_M, A_j)$ can be taken in its
cobordism class such that the
 characteristic class
$\kappa_M \in H^1(M^{n-k};\Z/2)$ of the skew-framing $\Xi_M$ is
given by the following composition:
$$\kappa_M : M^{n-k} \to \RP^{n-k-i} \subset \RP^{\infty},$$
$i < n-k$. We shall denote the map $M^{n-k} \to \RP^{n-k-i}$
described above again by $\kappa_M$.

 Let us consider the manifold
 $Q^i \subset M^{n-k}$, given by the formula:
\begin{eqnarray}\label{Q}
Q^i = \kappa_M^{-1}(pt), \quad pt \in \RP^{n-k-i}.
\end{eqnarray}
The manifold is equipped with the natural framing $\Psi_Q$,
because the restriction of the skew-framing $\Xi_M$ over the
submanifold $Q^i \subset M^{n-k}$ is a framing, which is denoted by $\Psi_Q$.

Moreover, the restriction of cohomology classes of the collection
$A_j(M)$, restricted over the submanifold $Q^i \subset M^{n-k}$, determines the
collection $A_j(Q)$ of cohomology classes on $Q^i$. Note that the
class $\kappa_0$ in the collection $A_j(Q)$ is the trivial class. The
immersion $\varphi_Q: Q^i \looparrowright \R^n$ is defined as the
restriction of the immersion $\varphi$ over the submanifold $Q^i
\subset M^{n-k}$. A triple $(\varphi_Q, \Psi_Q, A_j(Q))$
determines an element $J_{sf;A_j}^i(x)=y \in
Imm^{sf;A_j}(i,n-i)$.

\begin{lemma}\label{compspes}

The element  $x=[(\varphi, \Xi_M, A_M)] \in Imm^{sf,A_j} (n-k,k)$, which
admits a compression of the order $(i-1)$, admits  a
compression of the order $i$ if and only if the element
$J_{sf;A_j}^i(x)=y =[(\varphi_Q, \Xi_Q, A_j(Q))] \in
Imm^{sf,A_j}(i,n-i)$ is trivial.
\end{lemma}

\subsubsection*{Proof of Lemma $\ref{compspes}$}

At the first step let us prove that if $y=0$ then $x$ admits a
compression of the order $i$. Let us consider a skew-framed in the
codimension $(n-i)$ $(i+1)$--dimensional manifold
$(P^{i+1},\Xi_P)$ with boundary $\partial P^{i+1} = Q^i$, equipped
with the collection $A_j(P)$ of cohomology classes, such that the
restriction of the skew-framing $\Xi_P$ over the boundary $Q^i$ is
a framing coincided with the framing $\Psi_Q$ and the restriction
of the collection $A_j(P)$ over the boundary $Q^i$ coincides with
the collection $A_j(Q)$.

Let us describe a normal surgery of the skew-framed manifold
$(M^{n-k},\Xi_M)$ into a skew-framed manifold $(T^{n-k},\Xi_T)$.
Let us construct a manifold with boundary called the body of a
handle. Let us consider the manifold $P^{i+1}$ and let us denote
the $(n-i)$--dimensional normal bundle over $P^{i+1}$  by $\nu_P$.
The normal bundle $\nu_P$ is equipped with the skew-framing
$\Xi_P$, i.e. the bundle map (an isomorphism on each fibre)
\begin{eqnarray}\label{nuP}
\nu_P \to (n-i)\kappa_P
\end{eqnarray}
is well-defined. Let us denote by $U_P$ the disk bundle over $P$
spanned by the first $(n-k-i)$  factor in the Whitney sum
$(\ref{nuP})$. The manifold $U_P$ is a skew-framed manifold in the
codimension $k$  with boundary, this manifold  will be called the body
of a handle.

The boundary $\partial U_P$ of the body of the handle contains a
submanifold $Q^i \times D^{n-k-i}$,  the total space of the disk
bundle over the manifold $Q^i$. Let us consider the Cartesian
product $M^{n-k} \times I$ of the manifold $M^{n-k}$ and the unite
segment $I=[0,1]$.
 A second copy of the manifold
 $Q^i \times D^{n-k-i}$ is embedded into the submanifold  $M^{n-k} \times \{1\} \subset \partial(M^{n-k}
\times I)$, this is a regular neighbourhood of the submanifold $Q^i
\times \{1\} \subset M^{n-k} \times \{1\}$. Let us define the
manifold  $T^{n-k}$ by the following formula:
\begin{eqnarray}\label{T}
 T^{n-k}= \partial^+((M^{n-k} \times I) \cup_{Q^i \times D^{n-k-i}}  U_P),
\end{eqnarray}
where by $\partial^+$ is denoted the "upper" component of the
cobordism i.e. the component that contains the last part $\partial
U_P \setminus (Q^i \times D^{n-k-i})$  of the boundary of the body
$U_P$.

After the standard operation called "smoothing the corners" the
$PL$--manifold  $T^{n-k}$ becomes a smooth closed smooth manifold.
The immersion $\varphi_T: T^{n-k} \looparrowright \R^n$ (this
immersion is well-defined up to a regular homotopy), the skew
framing  $\Xi_T$  with the characteristic class $\kappa_T$ (i.e.
the bundle fibrewised isomorphism $\nu_T \to k\kappa_T$) are
well-defined. The manifold $T^{n-k}$ is equipped with the
collection $A_j(T)$ of the collection of characteristic classes,
each class in the collection is determined by the gluing of the
corresponded classes of the two components in the decomposition
$(\ref{T})$. The class $\kappa_0(T)$ of the collection $A_j(T)$
coincides with the characteristic class $\kappa_T$ of the
skew-framing $\Xi_T$. The triple $(\varphi_T, \Xi_T, A_j(T))$
determines an element in the cobordism group $
Imm^{sf;A_j}(n-k,k)$ and by the construction $[(\varphi, \Xi_M,
A_j)] = [(\varphi_T, \Xi_T, A_j(T))]$.

Let us prove that the element $[(\varphi_T, \Xi_T, A_j(T))]$
admits a compression of the order   $i$. Let us prove that the
characteristic class $\kappa_T$ is represented by a classifying
map $\kappa_T: T^{n-k} \to \RP^{n-k-i-1} \subset \RP^{\infty}$.
Take a positive integer $b$ big enough and let us consider
$\kappa_M : M^{n-k} \to \RP^b \subset \RP^{\infty}$, such that
$$\kappa_M^{-1}(\RP^{b-n+k+i}) = Q^i ,$$
$Q^i \subset M^{n-k}$, $\RP^{b-n+k+i} \subset \RP^b$.

Let us consider the mapping $g: P^{i+1} \to \RP^{b-n+k+i}$, the
restriction of this mapping to the component of the  boundary
$\partial P^{i+1} = Q^i$ coincides with the map $\kappa_M
\vert_{Q^i}$ (the mapping to a point).  Let us consider the "thickening"  $h: U_P \to
\RP^b$ of the map $g$, this map $h$ is defined by the standard
extension of the map $g$ to the body of the handle $U_P$ using the standard coordinate system in a regular neighbourhood of $\RP^{b-n+k+i} \subset \RP^b$.

The map  $g': M^{n-k} \cup_{Q^i \times D^{n-k-i}} U_P \to \RP^b$,
$g' \vert_{U_P} = g$ is well-defined and the restriction
 $g' \vert_{T^{n-k} \subset
M^{n-k} \cup U_P}$ does not meet the submanifold
 $\RP^{b-n-k+i} \subset \RP^b$.
 The space
$\RP^{b} \setminus \RP^{b-n-k+i}$ is retracted to the its subspace
$\RP^{n-k-i-1}$ by a deformation, the required compression of the
map $\kappa_T$ of the order  $i$ is constructed. We have proved
that the element $x$ admits a compression of the order $i$.

Let us prove the inverse statement: assume that the element
$x=[(\varphi, \Xi_M, A_j)]$  admits a compression of the order
$i$, then the triple $(\varphi_Q,\Psi_Q, A_j(Q)))$, $Q^i$ is given
by the equation $(\ref{Q})$ determines the trivial element in the
cobordism group $Imm^{sf;A_j}(i,n-i)$.

Let  $(\varphi_W,\Xi_W,A_j(W))$ be a triple, where $W^{n-k+1}$ is
a manifold with boundary, $\partial W^{n-k+1} = M^{n-k} \cup
M_1^{n-k}$; $(\varphi_W, \Xi_W)$ is a skew-framed immersion of the
manifold $W$ into $\R^n \times I$; $A_j(W)$ is a collection of
characteristic classes. Moreover, the triple
$(\varphi_W,\Xi_W,A_j(W))$ determines a cobordism between the
triples $(\varphi_M, \Xi_M, A_j(M))$ and $(\varphi_{M'}, \Xi_{M'},
A_j(M'))$, where the pair $(M'^{n-k}, \kappa_{M'})$ (the
cohomology class $\kappa_{M'}$ is the characteristic class of the
skew-framing $\Xi_{M'}$ and this class is included into the
collection $A_j(M')$) admits a compression of the order $i$, i.e.
the classifying map $\kappa_{M'} = \kappa_W \vert_{M'}$ is given
by the following composition:
$$ \kappa_{M'}: M'^{n-k} \to \RP^{n-k-i} \subset \RP^{\infty}.$$

Let us consider the standard submanifold  $\RP^{b-n-k+i} \subset
\RP^b$, this submanifold intersects the submanifold  $\RP^{n-k-i}
\subset \RP^b$ at a point $pt \in \RP^{n-k-i}\setminus
\RP^{n-k-i-1}$  and does not intersect the standard  submanifold
$\RP^{n-k-i-1} \subset \RP^{n-k-i}$. The image
$Im(\kappa_M(M^{n-k}))$ is in the submanifold $\RP^{n-k-i} \subset
\RP^b$, the image  $Im(\kappa_{M'}(M'^{n-k}))$ is in the
submanifold $\RP^{n-k-i-1} \subset \RP^{n-k-i}$.

Let us denote by  $P^{i+1}$ the submanifold
$F^{-1}(\RP^{b-n+k+i})$ (we assume that $F$ is transversal along
the submanifold $\RP^{b-n+k+i} \subset \RP^{b}$). By the
construction $\partial P^{i+1} = Q^i$. Let us define a
skew-framing $\Xi_P$ in the codimension $(n-i)$ as the direct sum
of a skew-framing of the submanifold $P^{i+1} \subset W^{n-k+1}$
and the skew-framing $\Xi_W$, restricted to the submanifold
$P^{i+1} \subset W^{n-k+1}$.

The restriction of the skew-framing $\Xi_W$ on $\partial W^{n-k+1}
=Q^{n-k}$ coincides with the skew-framing $\Psi_Q$ with the
trivial characteristic class $\kappa_Q$ (i.e. the skew-framing
$\Psi_Q$ is the framing). The restriction of the collection
$A_j(P)$ of cohomology classes on $\partial W^{n-k+1} =Q^{n-k}$
coincides with the collection $A_j(Q)$. This proves that the
triple $(\varphi_Q,\Psi_Q, A_Q)$, is a boundary. Lemma
$\ref{compspes}$ is proved. \qed

\subsubsection*{Proof of Theorem $\ref{totalobstr}$}

Let us assume that a compression of the order $(i-1)$, $i<q$ for
an element $x \in Imm^{sf;A_j}(n-k,k)$, $x=[(\varphi, \Xi_M,
A_M)]$ is well-defined. By Lemma $\ref{compspes}$, the obstruction
to a compression of the order $i$ of the element $x \in
Imm^{sf;A_j}(n-k,k)$ represented by the same triple $(\varphi,
\Xi_M, A_M)$, coincides with the obstruction of a compression of
the same order $i$ for the element $J_{k'}^{sf}(x) \in
Imm^{sf;A_j}(n-k',k')$.  Therefore, by induction over $i$, the
total obstruction for a compression of the order $q$ for the
element $x$ is trivial if and only if the total obstruction for a
compression of the order $q$ for the element $J_{k'}^{sf}(x)$ is
trivial. Theorem $\ref{totalobstr}$ is proved.
\[  \]

To prove the Compression Theorem $\ref{comp}$ the following
construction by U.Koschorke of the total obstruction for a homotopy
of a bundle map into a bundle monomorphism on each fibre of the
bundles (see [K]) is required.

Let
 $\alpha \to Q^q$, $\beta \to Q^q$ be a pair of the vector bundles over the smooth manifold
 $Q^q$ (we do not assume that the manifold $Q^q$ is closed)
  $dim(\alpha)=a$,
$dim(\beta)=b$, $dim(Q^q)=q$,  $2(b-a+1)<q$. Let $u: \alpha \to \beta$ be a generic vector bundle morphism.
let us denote by
$\Sigma \subset Q^q$ a submanifold, given by the formula:
\begin{eqnarray}\label{Sigma}
\Sigma = \{ x \in Q^q \vert Ker(u_x: \alpha_x \to \beta_x) \ne 0
\}.
\end{eqnarray}
This manifold $\Sigma$ is the singular manifold of the bundle morphism $u$.
Note that under the presented dimensional restrictions, for a generic vector bundle morphism $u$
we have
$rk(u) \ge a-1$.The codimension of the submanifold
$\Sigma \subset Q^q$ is equal to
$b-a+1$.

Let us describe the normal bundle  of the submanifold
$(\ref{Sigma})$, this bundle will be denoted by $\nu_{\Sigma}$. Let us denote by $\lambda: E(\lambda) \to \Sigma$ the linear subbundle,
determined as the subbundle of kernels of the morphism $u$ over the singular submanifold
 $\Sigma \subset Q^q$. Therefore, the following inclusion of  bundles over $\Sigma$
  $\varepsilon: \lambda \subset \alpha$ is well-defined.
Let us denote by $\Lambda_{\alpha}$ the bundle over $\Sigma$, this
bundle is the orthogonal complement to the subbundle
$\varepsilon(\lambda) \subset \alpha$. A natural vector-bundle
morphism over $\Sigma$ (isomorphism of fibres)
 $v: \Lambda_{\alpha} \subset \beta$ is well-defined.
Let us define the bundle $\Lambda_{\beta}$ over $\Sigma$ as the
orthogonal complement to the subbundle $v(\Lambda_{\alpha})$ in
the bundle $\alpha \vert_{\Sigma}$. The normal bundle
$\nu(\Sigma)$ is determined by the following formula:
\begin{eqnarray}\label{016}
\nu(\Sigma) = \lambda \otimes \Lambda_{\beta}.
\end{eqnarray}

If the manifold
$Q^q$ has a boundary $\partial Q$ and the vector bundles morphism  $u$ is the
morphism of the bundles over the manifold with boundary, then the singular submanifold
$\partial \Sigma \subset \partial Q$ of the restriction $u \vert_{\partial Q}$
is a boundary of the submanifold $\Sigma \subset Q^q$ with the normal bundle, given by the same formula $(\ref{016})$.

In the paper [K] (in this paper there is a reference to the
previous papers by the same author) a cobordism group of
embeddings of of manifolds in $Q^q$ (in this construction the
manifold $Q^q$ is closed) of codimension
 $b-a+1$ with an additional structure of the normal bundle, given by the equation  ($\ref{016}$)
is defined. For an arbitrary generic vector bundle morphism $u:
\alpha \to \beta$ an element in this cobordism group is
well-defined. This element is the total obstruction of a homotopy
of the vector bundle morphism $u$ to a fibrewised monomorphism.

Let $\hat \nu: E(\hat \nu) \to \RP^{2^k-1}$ be a vector bundle, $\dim(\hat \nu)=
n+1-2^k$, $2^k < n+2$ over the standard projective space, isomorphic to the following Whitney sum:
\begin{equation}\label{hatnu}
\hat \nu_{\RP^{2^k-1}} = (n+1-2^k)\kappa_{\RP},
\end{equation}
 where  $\kappa_{\RP}$ is the canonical line bundle over
$\RP^{2^k-1}$. Let us denote the Whitney sum $\hat \nu_{\RP^{2^k-1}} \oplus \kappa_{\RP}$ by
$\nu_{\RP^{2^k-1}}$, we get:
\begin{eqnarray}\label{nu}
 \nu_{\RP^{2^k-1}}= (n-2^k+2)\kappa_{\RP}.
\end{eqnarray} 
  The standard projection
 \begin{eqnarray}\label{proj} 
   \pi: \nu_{\RP^{2^k-1}} \to
\hat \nu_{\RP^{2^k-1}}
\end{eqnarray}
 with the kernel $\kappa_{\RP}$ is well-defined.
 In the case
\begin{eqnarray}\label{b}
b(2^k) \le n+2
\end{eqnarray}
(the positive integer $b(r)$, $r=2^k$ equals to the minimal
power of $2$, such that the bundle $b(2^k) \kappa_{\RP^{2^k-1}}$ isomorphic to the trivial bundle; in the right hand side of (\ref{b}) we have a power of $2$), the bundle $\nu_{\RP^{2^k-1}}$  by (\ref{nu}) is isomorphic to the normal bundle
over the projective space $\RP^{2^k-1}$.

Let us give a definition.

\subsubsection*{Definition of an admissible family
of sections of the bundle $\hat \nu_{\RP}$}

We shall say that a generic $s$--family of sections
$$\hat \psi =
\{\hat \psi_1, \dots, \hat \psi_s\}, \quad \hat \psi: s\varepsilon \to \hat \nu_{\RP}$$
of the bundle $\hat \nu_{\RP}$ is admissible, if
there exists a regular $s$--family
$$\psi = \{
\psi_1, \dots, \psi_s\}, \quad \psi: s\varepsilon \to \nu_{\RP}$$
of sections of the bundle  $\nu_{\RP}$ satisfied the condition:
$\pi \circ \psi=\hat \psi$.
\[  \]

\begin{lemma}\label{addmiss}
Let us assume that $s \le n+2-2^{k+1}-k-1$ and $n \equiv -2
\pmod{2^{2^k}}$, $n >0$. Then the bundle $\hat \nu_{\RP^{2^k-1}}$ has an
admissible $s$-family of sections.
\end{lemma}

\subsubsection*{Proof of Lemma $\ref{addmiss}$}

By the Davis table  the projective space $\RP^{2^k-1}$ is
immersible into the Euclidean space $\R^{2^{k+1}-k-1}$ (this is
not the lowest possible dimension of the target Euclidean space of
immersions). By the equation $(\ref{b})$, and because $b(k) \le
2^k$ for $n=2^{2^k}-2 \pmod{2^{2^k+1}}$, the bundle
$\nu_{\RP^{2^k-1}}$, given by (\ref{nu}), admits a generic regular $s$--family of sections,
denoted by $\psi$. The projection $\hat \psi = \pi \circ \psi$ of
this regular family is the admissible $s$--family of sections of
the bundle $\hat \nu_{\RP}$. Lemma $\ref{addmiss}$ is proved. \qed

Let us consider an admissible generic $s$--family of sections $\hat \psi$ of the bundle $\hat \nu_{\RP}$.
Let us denote by $\Sigma \subset \RP^{2^k-1}$ the singular manifold of the family $\hat \psi$.
This denotation corresponds to $(\ref{Sigma})$, if we take $\alpha \equiv s \varepsilon$,
$\beta \equiv \hat \nu$, $u=\psi$. In the following lemma we will describe the normal bundle
$\nu_{\Sigma}$ of the submanifold $\Sigma \subset \RP^{2^k-1}$.

\begin{lemma}\label{7.5}
Let us assume that
\begin{eqnarray}\label{s}
s = n+2-2^{k+1}-k-1, \quad k \ge 2,
\end{eqnarray}
where $n = 2^{2^k}-2 \pmod{2^{2^k}}$. Then the submanifold $\Sigma
\subset \RP^{2^k-1}$ of  singular admissible  sections in a generic $s$--family of
the vector-bundle $\hat \nu_{\RP}$ $(\ref{hatnu})$ is a smooth submanifold of
dimension $k$, the normal bundle $\nu_{\Sigma}$ of the submanifold is equipped with a
skew-framing $\Xi_{\Sigma}: (2^k-1-k) \kappa_{\Sigma} \equiv
\nu_{\Sigma}$, the characteristic class  $\kappa_{\Sigma}=\kappa_{\Sigma}(\Xi)$ of this
skew-framing coincides with the restriction
$\kappa_{\RP}\vert_{\Sigma \subset \RP^{2^k-1}}$.
 \end{lemma}

\subsubsection*{Proof of Lemma $\ref{7.5}$}

 Let us describe the normal bundle $\nu_{\Sigma}$ by means of Koschorke's Theorem.
Moreover, let us define a skew-framing of this bundle.
Let us denote by
$\lambda \subset s\varepsilon $ the subbundle of the kernels of the family $\hat \psi$
over the submanifold $\Sigma$. By the assumption the the family $\hat \psi$ is admissible,
therefore:
\begin{eqnarray}\label{lambda}
\lambda = \kappa_{\Sigma}
\end{eqnarray}
(recall,  by $\kappa_{\Sigma}$ the vector bundle $\kappa_{\RP} \vert_{\Sigma}$ is denoted).
Let us denote the orthogonal complement to $\kappa_{\Sigma} \subset s \varepsilon$ over $\Sigma$ by  $s\varepsilon -\kappa_{\Sigma}$.

By the construction:
\begin{eqnarray}\label{barPsi}
(s\varepsilon - \kappa_{\Sigma}) \oplus \kappa_{\Sigma} = s \varepsilon.
\end{eqnarray}
Let us prove that the bundle  $\nu_{\Sigma}$ satisfies the equation:
\begin{eqnarray}\label{lambda1}
 \nu_{\Sigma} \equiv (n+2 - 2^k -k-1) \kappa_{\Sigma}.
 \end{eqnarray}

Let us denote the orthogonal complement to $\hat \psi(s\varepsilon -\kappa_{\Sigma})$ in the bundle
$\hat \nu_{\RP^{2^k-1}}\vert_{\Sigma}$ by $\Lambda$. Obviously, $dim(\Lambda)=2^k-1+k$.
 By the Koschorke Theorem the following isomorphism of vector bundles is well-defined:
 $$\Lambda \otimes \kappa_{\Sigma} \equiv \nu_{\Sigma}.$$
The bundle $\Lambda$ is a stable bundle over $\Sigma$, which is stably trivial bundle with the prescribed stable trivialization. This is followed from the fact: the stable vector bundle in the left hand side of the following equation is a  trivial bundle (with a prescribed stable trivialization): 
$$\hat{\nu}_{\RP^{2^k-1}} \oplus \kappa_{\Sigma} = \hat{\psi}(s\varepsilon - \kappa_{\Sigma}) \oplus \kappa_{\Sigma} \oplus \Lambda.$$  
 This proves the equation $(\ref{lambda1})$.
   The isomorphism $(\ref{lambda1})$ defines a skew-framing $\Xi_{\Sigma}$
of the bundle $\nu_{\Sigma}$   with the characteristic class $\kappa_{\Sigma}$.
The Lemma $\ref{7.5}$ is proved. \qed

 \subsubsection*{Remark}
Because the restriction of the normal bundle  $\nu_{\RP}$ over
$\Sigma$ is isomorphic to the trivial bundle, using the following stable
isomorphism over $k$-dimensional manifold $\Sigma$:
$$ b(k)\kappa_{\RP} = b(k)\varepsilon,$$ 
the skew-framing $\Xi_{\Sigma}$ determines the
skew-framing of the normal bundle (in the Euclidean space) of the
manifold $\Sigma$.
\[  \]

Let us consider an element
$$ x \in Imm^{sf,A_j}(2^{k}-2,n-2^{k}+2), \quad n=-2(mod 2^{k+1}),$$
given by the cobordism class of a triple $(\varphi, \Xi_M, A_j(M))$. Put $m=2^{k}-2$.
Let us consider the vector bundle
$\nu_{\RP^{2^k-1}} \to \RP^{2^k-1}$, $\dim(\nu)=n-m$.
The normal bundle
$\nu_M$ (\ref{nu}) of the manifold $M^{m}$ is given by the formula:
$$\nu_M = \kappa_M^{\ast}(\nu_{\RP^{2^k-1}}).$$
Let us consider the map
\begin{eqnarray}\label{lam}
\lambda_M: M^{m} \to \RP^{2^k-1} \times \prod_{i=1}^j \RP^{\infty},
\end{eqnarray}
see $(\ref{lambdaM})$, constructed by means of the collection
$A_j$. Let us consider the standard projection  $\pi_0:
\RP^{2^k-1}(0) \times \prod_{i=1}^j \RP^{\infty}(0) \to \RP^{2^k-1}$ on
the factor $\RP^{2^k-1}(0)$. The composition  $\pi_0 \circ \lambda_M:
M^m \to \RP^{2^k-1}(0)$ coincides with the map $\kappa_M$.

Let us define the subbundle $\hat \nu_M \subset \nu_M$ of the
codimension 1 (i.e. of the dimension $\dim(\hat \nu_M)=n+1-2^{k}$)
by the formula:
$$\hat \nu_M = \kappa_M^{\ast}(\hat \nu_{\RP^{2^k-1}}).$$
Let us define a family of sections  $\hat \xi_M=\{\hat \xi_1,
\dots, \hat \xi_s\}$, $s=n+3-2^{k+1}+k$ of the bundle $\hat
\nu_M$. This collection is the pull-back image of an admissible
collection of sections $\hat \psi = \{\hat \psi_1, \dots \hat
\psi_s\}$ of the bundle $\hat \nu_M$ by the map $\kappa_M$. Let us
denote by
\begin{eqnarray}\label{NsubM}
N^{k-1} \subset M^m
\end{eqnarray}
the submanifold  of singular sections. It is not hard to prove
that $\dim(N^{k-1})=(k-1)$. The manifold $N^{k-1}$ is equipped
with the collection $A_j(N)$ of cohomology classes from the group $H^1(N^{k-1};\Z/2)$, a class of
$A_j(N)$ is defined as the restriction of the corresponded class
of the collection $A_M$ over the submanifold  $N^{k-1} \subset
M^m$. The class $\kappa_M \vert_N$ in the collection $A_j(N)$ is
denoted by $\kappa_N$. Let us denote the immersion $\varphi
\vert_N$ by $\varphi_N$.

Let us denote by $\nu_N$ the normal bundle of the immersed
(embedded by general position arguments) manifold
$\varphi_N(N^{k-1})$ in the Euclidean space $\R^n$. The normal
bundle $\nu_N$ is isomorphic to the Whitney sum  $\nu_N = \nu_M
\vert_N \oplus \nu_{N \subset M}$, where by $\nu_{N \subset M}$ is
denoted the normal bundle of the submanifold $N^{k-1} \subset M^m$
inside the manifold $M^m$.

By Lemma $\ref{7.5}$ and by the transversality of the map
$\kappa_M$ along the submanifold $\Sigma \subset \RP^{2^k-1}$, the
bundle $\nu_{N \subset M}$ is equipped with the skew-framing
$\Xi_{N \subset M}$ with the characteristic class $\kappa_N$. The
bundle $\nu_M \vert_N$ is also equipped with a skew-framing by the Koschorke Theorem.
The calculation of the characteristic class  of this skew-framing is analogous to Lemma \ref{7.5}. This gives a skew-framing  $\Xi_N$ of the
immersion $\varphi_N$ of codimension $(n-k+1)$.

\begin{lemma}\label{7.6}
The triple $(\varphi_N, \Xi_{N}, A_j(N))$ determines an element
$x_{k-1} \in Imm^{sf,A_j}(k-1,n-k+1)$, this element is the total
obstruction of a compression of the order $(k-1)$ for the element
$x \in Imm^{sf,A_j}(2^{k}-2,n-2^{k}+2)$.
\end{lemma}

\subsubsection*{Proof of Lemma $\ref{7.6}$}

Let us consider the submanifold $\Sigma^k \subset \RP^{2^{k}-1}$
of singularities of the family of sections  $\hat \psi$, let us
re-denote this manifold by $\Sigma_0^{k}$. This manifold is
equipped by the following natural stratification (a filtration):
\begin{eqnarray}\label{017}
\emptyset \subset \Sigma_{k}^0 \subset \dots \subset
\Sigma_{1}^{k-1} \subset \Sigma^{k}_0 \subset \RP^{2^{k}-1}.
\end{eqnarray}

The submanifold $\Sigma_i$, $dim(\Sigma_i)=k-i$ in $(\ref{017})$
is defined as the singular submanifold of the subfamily of the
first $(s-i)$ sections in $\hat \psi$. A straightforward
calculations implies that the fundamental class
 $[\Sigma_i]$ of a submanifold in (\ref{017}) represents in the group
$H_{k-i}(\RP^{2^{k}-1};\Z/2)$  the only generator of
this group: this homology class is dual to the characteristic
class $\bar w_{n+1-2^{k}-k+i}[(n+1-2^{k})\kappa_{\RP}]$ of the
bundle $\hat \nu_{\RP^{2^k-1}}$.

Without loss of a generality we assume that the map
$\kappa_M : M^m \to \RP^{2^{k}-1}$ is transversal along
the stratification  $(\ref{017})$.
Let us denote the inverse image of the stratification
 $(\ref{017})$ by
\begin{eqnarray}\label{018}
 N_{k-1}^0 \subset N_{k-2}^1 \subset \dots \subset N_{0}^{k-1}
\subset M^m.
\end{eqnarray}
The top manifold $N_0^{k-1}$ of the filtration  $(\ref{018})$  coincides with manifold
$N^{k-1}$, defined above.

Let us prove the lemma by the induction over the parameter $i$,
$i=0, \dots, k-1$. Let us assume that the image of the map
$\kappa_M : M^m \to \RP^{2^{k}-1}$ is in the standard projective
subspace $\RP^{2^{k}-2-i} \subset \RP^{2^{k}-1}$.  In this case
$N^{i-1}_{k-i}=\emptyset$. By the standard argument we may assume
that the stratum $\Sigma_{k-i}^i$  intersects in the general
position the standard submanifold $\RP^{2^{k}-1-i}$ of the
complementary dimension at the only point. (The index of
self-intersection of this two submanifolds in the manifold
$\RP^{2^{k}-1}$ is odd and well-defined modulo 2.)

The framed manifold  $N_{k-i-1}^i$ is the regular preimage of the
marked point by the map $\kappa_M$ (the image of this map is in
the submanifold
 $\RP^{2^{k}-2-i}$). Let us denote an element represented by the cobordism class of the triple $(\varphi_{N_{k-i-1}},\Xi_{N_{k-i-1}},
A_j({N_{k-i-1}}))$ in the cobordism group $Imm^{sf;A_j}(i,n-i)$ by $x_i$. By Lemma $\ref{compspes}$ the condition
 $x_i=0$, $i \le k-2$ is satisfied if and only if there exists a normal cobordism
of the map $\kappa_M$ to the map $\kappa_{M'}: M'^{m} \to
\RP^{2^{k}-2-i} \subset \RP^{2^{k}}$.
Therefore a
compression of the order $i+1$, $i+1 \le k-1$ of an element $x$ is
well defined. If we put $i+1 = k-1$ we have a compression of the
order $(k-1)$. Lemma $\ref{7.6}$ is proved. \qed
 \[  \]

The following proposition is the main step in the proof of Theorem $\ref{comp}$.

\begin{proposition}\label{38}
	
Let $n=-2(mod(2^k))$, $n>2^k$, $x \in Imm^{sf;A_j}(2^k-2,n-2^k+2)$ be an arbitrary element in the kernel of the homomorphism 
$(\ref{Omega})$ (in this formula we put $b(k)=2^k$, but this is not a minimal value for $b(k)$, see $(\ref{b})$). Let $x_{k-1} \in Imm^{sf;A_j}(k-1,n-k+1)$
be the total obstruction for a compression of the order $(k-1)$ of the element
$x$. Then the element
$x_{k-1}$ is in the image of the transfer homomorphism $(\ref{transfer})$, i.e. there exists an element $y_{k-1} \in
Imm^{sf;A_{j+1}}(k-1,n-k+1)$, for which
$r^!_{j+1}(y_{k-1})=x_{k-1}$. 
\end{proposition}

\subsubsection*{Proof of Proposition $\ref{38}$}

Let us assume that the cobordism class of the element $x$ is given
by a triple $(\varphi_M,\Xi_M,A_j(M))$, where $\varphi_M: M^{m}
\looparrowright \R^n$ is an immersion, $\dim(M^{m})=m=2^k-2$. 
Because $codim(\varphi_M)$ is odd, $M^m$ is oriented. 

Let
us consider the normal bundle $\nu_M \cong (n-2^{k}+2)\kappa_M$
over the manifold $M^{m}$ and the subbundle  $\hat \nu_M \subset
\nu_M$ of the codimension 1, $\hat \nu_M \equiv
(n-2^{k}+1)\kappa_M$. Let us prove that there exists a regular
$s$--family of sections $\hat \psi$ of the bundle $\hat \nu_M$,
$s=n+2+k - 2^{k+1}$.

Let us denote by $D(\kappa_M)$ a manifold with boundary, the total
space of the disk bundle, associated with line bundle $\kappa_M$.
The vector bundle $\hat \nu_M$ is lifted to the vector bundle over
$D(\kappa_M)$ (we will denote this lift again by  $\hat \nu_M$).
By the R.Cohen theorem [C] there exists an immersion  $D(\kappa_M)
\looparrowright \R^{2^{k+1}-2-k}$, because
$\dim(D(\kappa_M))=2^k-1$ and $\alpha(2^k-1)=k$.
An elementary self-contained proof of the statement is
below in Theorem \ref{39}.

Equivalently, the
bundle $\hat \nu_M$ admits an $s$--family of regular sections. The
tautological lift, denoted by $\psi$, of the regular $s$--family
$\hat \xi$ of the bundle $\hat \nu_M$ to a regular $s$--family of
sections of the bundle $\nu_M$ is defined.

Let us consider the admissible $s$--family of sections $\hat \psi$
of the bundle $\hat \nu_M$. This family of sections is defined as
the pull-back of admissible $s$--family of sections of the bundle
$\hat \nu_{\RP^{2^k-1}}$, see Lemma $\ref{addmiss}$. A regular $s$--family
of sections $\psi$ of the bundle $\nu_M$ is defined as the lift of
the admissible $s$--family $\hat \psi$.

Let us consider the manifold $M^m \times I$ and let us define the
bundle $\nu_{M \times I}$ by the formula $\nu_{M \times I} = p_M^{\ast}(\nu_M)$, where $p_M: M^m \times I \to M^m$
is the standard projection on the second factor.
Let us define the
bundle $\hat \nu_{M \times I}$ by the formula $\hat \nu_{M \times I} = p_M^{\ast}(\hat \nu_M)$.
Let us consider a generic $s$--family of sections
$X= \{\chi_1, \dots \chi_s\}$ of the bundle  $\nu_{M \times I}$
with the following boundaries conditions:
\begin{eqnarray}\label{220}
X = \psi \quad over \quad M^{m} \times \{1\},
\end{eqnarray}
\begin{eqnarray}\label{221}
X = \xi \quad over \quad M^{m} \times \{0\}.
\end{eqnarray}

Let us denote by $V^{k-1} \subset M^{m} \times I$ the singular subset of the family $X$.
By the general position argument this subset is a closed submanifold in $M^{m} \times I$,
because over the boundary $\partial(M^{m} \times I)$ the family $X$ is regular.
Let us denote by
$$\hat X \{\hat \chi_1, \dots, \hat \chi_s \}$$
the projection of the $s$--family $X$ into the $s$--family of sections of the vector bundle
$\hat \nu_{M \times I}$. The $s$--family $\hat X$ satisfies the following boundary conditions:
\begin{eqnarray}\label{222}
\hat X = \hat \psi \quad over \quad M^{m} \times \{1\},
\end{eqnarray}
\begin{eqnarray}\label{223}
\hat X = \hat \xi \quad over \quad M^{m} \times \{0\}.
\end{eqnarray}

Let us denote by $\hat K^{k} \subset M^{m} \times I$ the subset of
singular sections of the $s$--family $\hat X$. This subset $\hat
K^k$ is a $k$--dimensional manifold with boundary. The only
component of the boundary $\partial \hat K^k$ is a submanifold of
$M^{m} \times \{1\})$, this component coincides with the
submanifold $N^{k-1} \subset M^m \times \{1\}$ of the
singularities of the admissible family $\psi$, see
$(\ref{NsubM})$ (recall, the family  $\hat{\xi}$ is regular and  is admissible).
By definition the following inclusion is well-defined: $V^{k-1} \subset \hat{K}^k$.

By Lemma $\ref{7.5}$, the triple $(\varphi_N, \Xi_N, A_N)$ is
well-defined. Here $\varphi_N = \varphi \vert_N$, $\Xi_N$ is the
skew-framing of this immersion, $A_j(N)$ is the restriction of the
collection  $A_j(M)$ on $N^{k-1} \subset M^m \times \{1\}$.
 The triple $(\varphi_N, \Xi_N, A_j(N))$ represents the total obstruction $x_{k-1} \in Imm^{sf;A_j}(k-1,n-k+1)$
for a compression of the order  $(k-1)$ of the element
$x=[(\varphi, \Xi_M, A_j(M))]$, $x \in Imm^{sf;A_j}(2^k-2,n-2^k+2)$.

Let us use the formula
$(\ref{016})$ to calculate the normal bundle of the submanifold
$\hat K^{k} \subset M^m \times I$ and of the line normal bundle of the submanifold
$V^{k-1} \subset \hat K^{k}$.

Let us denote by $\lambda$ the line bundle over $V^{k-1}$ of kernels of the $s$--family of sections $X$.
Let us prove that the normal bundle $\nu_{V}$ of the submanifold
$V^{k-1} \subset M^m \times I$
 is given by the formula:
\begin{eqnarray}\label{normV}
\nu_{V} \equiv \varepsilon \oplus  (m-k+1)\lambda, \quad m=2^k-2.
\end{eqnarray}
The restriction of the normal bundle $\nu_{M \times I} \vert_{V}$
over the submanifold $V^{k-1}$ is isomorphic to the bundle
$(n+2-2^{k})\kappa_{M \times I}$, and therefore, because $b(k-1)
\le 2^k$, is isomorphic to the trivial bundle:
\begin{eqnarray}\label{norm}
(n+2-2^{k})\kappa_{M \times I} \equiv (n+2-2^{k})\varepsilon.
\end{eqnarray}
This isomorphism is canonical, i.e. does not depends of $M^m$ and
of $V^{k-1}$. Let us define the vector bundle $\Lambda$ as the
orthogonal complement to the line subbundle $\lambda$ in the
trivial bundle of linear combinations of the base sections. The
bundle $\Lambda$ represents the stable vector bundle $-
\lambda$. Therefore the orthogonal complement of the subbundle
$\Xi(\Lambda)$ in the vector bundle $\nu_{M \times I}$ is the
subbundle in the vector bundle $\nu_M \vert_V$ isomorphic to the
vector bundle $\lambda \oplus (2^{k} - k -1)\varepsilon$. By the
Koschorke Theorem, the normal bundle of the submanifold $V^{k-1}
\subset M^{2^{k}-2} \times I$ is isomorphic to the vector bundle:
\begin{eqnarray}\label{normV}
\nu_{V \subset M \times I} \equiv \lambda \otimes(\lambda \oplus (2^{k}-k-1)\varepsilon)
\equiv \varepsilon \oplus (2^{k}-k-1)\lambda.
\end{eqnarray}
 This vector bundle
$\nu_{V \subset M \times I}$ also represents the stable normal
bundle of the manifold $V^{k-1}$ because of the equation
$(\ref{norm})$. The formula $(\ref{normV})$ is proved.

The restriction of the immersion $\varphi \times Id\vert_V$ is
 an immersion (an embedding)
 $\varphi_V: V^{k-1} \looparrowright \R^n \times [0,1]$. This immersion is a skew-framed immersion in the target spece $\R^n$, by the formula (\ref{normV}). By the computation the normal bundle of the immersion
 $\varphi_V$ is equipped with a skew-framing, denoted by $\Xi_V$. The collection of cohomology classes $A_j(V)$
 is defined by the formula $A_j(V) = A_j(M \times I) \vert_V$, where $A_j(M \times I)$ is induced from the given collection $A_j(M)$ of cohomology classes on $M^m$ by the projection $p_M$.
Let us define the collection of cohomology classes
$A_{j+1}(V)$, adding to the collection $A_j(V)$ the last cohomology class  $\kappa_{j+1} = \lambda \otimes \kappa_{M \times I}$.
The triple
 $y_{k-1}=(\varphi_V, \Xi_V, A_{j+1}(V))$ determines an element in the cobordism group $Imm^{sf;A_{j+1}}(k-1,n-k+1)$.

Let us denote by $U_V \subset \hat K^{k}$ a small closed regular
neighbourhood of the submanifold  $V^{k-1} \subset \hat K^{k}$. The
line bundle of kernels of the family $\hat X$ of sections over the
submanifold (with boundary) $U_V \subset M^m \times I$ is isomorphic
to the line bundle $p_{U_V,V}^{\ast}(\lambda)$, where $p_{U_V,V}:
U_V \to V$ is the projection of the neighbourhood on the central
submanifold. Let us denote the line bundle
$p_{U_V,V}^{\ast}(\lambda)$  by $\hat \lambda$. The orthogonal
complement to the subbundle $\hat \lambda$ in the vector bundle of
the linear combinations of the base sections over $U_V$ will be
denote again by $\hat \Lambda$.

Let us consider the subbundle $\hat X(\hat \Lambda)$ in the vector
bundle $\hat \nu_{M \times I}$, by this bundle we means the image of $\Lambda$ by the family of sections. Analogous calculations,
using the canonical isomorphisms of the vector bundles over $U_V$:
$(2^{k-1})\lambda \equiv (2^{k-1})\kappa_M \equiv
(2^{k-1})\varepsilon$, shows that the orthogonal complement in $\hat \nu_{M
\times I}$ of the subbundle $\hat X(\hat \Lambda)$ (this complement is denoted by
$-\hat X(\hat \Lambda)$), is isomorphic to the bundle $\lambda \oplus
(2^{k-1}-k-1)\varepsilon \oplus (2^{k-1}-1)\kappa_M$.  By the
Koschorke Theorem, the stable isomorphism class of the normal
bundle $\nu_{U_V \subset M \times I}$ of the submanifold $U_V
\subset M^{m} \times I$ is given by the formula:

$$\nu_{U_V \subset M \times I} \equiv (-k-1)\hat \lambda \oplus (-\kappa_M \otimes
\hat \lambda).$$

In particular, using (\ref{normV}),  the line normal bundle of the submanifold
$V^{k-1} \subset K^{k}$ is isomorphic to the line bundle  $\lambda \otimes \kappa_M$.

Let us denote $\partial U_V$ by $Q^{k-1}$. The space $Q^{k-1}$ is a closed manifold, $\dim(Q^{k-1})=k-1$.
The normal bundle $\nu_{Q \subset M\times I}$ of the submanifold $Q^{k-1} \subset M^m \times I$ is given
by the formula:
$$ \nu_{Q \subset M\times I} \equiv \varepsilon \oplus (m-k+1) \kappa_M. $$
The restriction of the immersion $\varphi \times id: M^m \times I \looparrowright \R^n \times I$
to the submanifold $Q^{k-1} \subset M^m \times I$ is regular homotopic to a skew-framed immersion
$\hat \varphi_Q: Q^{k-1} \looparrowright \R^n \times \{1\}$ оf the codimension $(n-k+1)$
with the skew-framing, denoted by $\hat \Xi_Q$, and with the characteristic class of this skew-framing $\hat \kappa_Q = \kappa_{M \times I} \vert_Q$. The manifold $Q^{k-1}$ is equipped by the collection $\hat A_j(Q)$ of cohomology classes, $\hat A_j(Q) = A_j(M \times I) \vert_{Q \subset M \times I}$. The triple
 $(\varphi_Q, \Xi_Q, A_j(Q))$ determines an element in the cobordism group  $Imm^{sf;A_j}(k-1,n-k+1)$.

The manifold $K^{k} \setminus U_V$ has a boundary consists of
the two components:
 $\partial(K^{k} \setminus U_V)= Q^{k-1} \cup K^{k-1}$.
 The restriction of the immersion $\varphi \times Id \vert_{(K^{k} \setminus U_V)}$ is regular homotopic to an immersion
 $\varphi_K: K^k \looparrowright \R^n \times I$ with the following the boundary conditions:
 $\varphi_K \vert Q = \varphi_Q: Q^{k-1} \looparrowright \R^n \times \{1\}$,
 $\varphi_K \vert N = \varphi_N: N^{n-1} \looparrowright \R^n \times \{0\}$.
 The immersion $\varphi_K$ is a skew-framed immersion with a skew-framing $\Xi_K$ and with the characteristic class $\kappa_K = \kappa_{M \times I} \vert_K$ of this skew-framing. The manifold $K^{k}$ is equipped by the collection of cohomology classes $A_j(K) = A_j(M \times I) \vert_{K \subset M \times I}$.

The triple
 $x_{k-1}=(\varphi_N, \Xi_N, A_j(N))$ determines an element in the cobordism group  $Imm^{sf;A_j}(k-1,n-k+1)$.
 The element $x_{k-1}$ is the total obstruction to a compression of the element $x=[(\varphi,\Xi_M,A_j)]$ of the order $(k-1)$.
The element $[(\varphi_Q, \hat \Xi_Q, \hat A_j(Q))]$ is the image by the transfer homomorphism $r_{j+1}^!$ of the element
 $y_{k-1}=[(\varphi_V, \Xi_V, A_{j+1}(V))]$. The elements   $x_{k-1}=[(\varphi_N, \Xi_N, A_j(N))]$ and $r^!_{j+1}(y_{k-1})=[(\varphi_Q, \hat \Xi_Q, \hat A_j(Q))]$ are equal. The cobordism between the elements $x_{k-1}$ and $r^!_{j+1}(y_{k+1})$ is given by the triple $(\varphi_K, \Xi_K, A_j(K))$.

Assume that $k-1$ is even. Let us prove that $y_{k-1}$ is in the kernel of $(\ref{Omega})$ (in this formula we assume that $b(k)=k+1$). The construction of  $y_{k-1}$ 
from $x$ is generalized for the image of $x$ by $(\ref{Omega})$ (in this formula $b(k)=2^k$). Because $x$ is in the kernel of $(\ref{Omega})$,
$y_{k-1}$ is also in the kernel of $(\ref{Omega})$.
Proposition $\ref{38}$ is proved. \qed

\subsubsection*{Proof of the Compression Theorem $\ref{comp}$}

Let us define a positive integer $\psi = \psi(q)$, by the formula
 $(\ref{psi})$ for $n-k=q-1$. By Proposition
$\ref{totaltrans}$ the total transfer homomorphism $(\ref{104})$, which is defined
on the group $Imm^{sf;A_{\psi}}(d-1,n-d+1)$ is the trivial
homomorphism. Let us define a positive integer  
\begin{eqnarray}\label{l}
l(d) =
\exp_2(\exp_2 \dots \exp_2(d-1+1)\dots +1), 
\end{eqnarray}
where the number of the
iterations of the function $\exp_2(x+1) = 2^{x+1}$ equals to
$\psi$, defined in Proposition \ref{totaltrans}, and the base of the iteration sequence is $x=d-1$.

Let  $l'$ be an arbitrary power of $2$, $l' \ge l(d)$. Let us
define  $n=l'-2$. Let us prove that an arbitrary element in
$Imm^{sf}(n-1,1)$ admits a compression of the order $d-1$.

Let us define $n_0=l(d)-2$, by the assumption $n_0 \le n$. Let us
consider the following sequence of $\psi$ integers: $2n_1 =
log_2(n_0+2)-2$, $2n_2= log_2(n_1+2)-2$, $\dots,$
 $2n_{\psi}=log_2(n_{\psi-1} +2)-2$. All this integers are positive and $n_{\psi}=d-1$.

Let $J^{sf}_{n_0}(x) \in Imm^{sf}(n_0,n-n_0)$ be the image an  element
$x \in Imm^{sf}(n-1,1)$ by $J^{sf}: Imm^{sf}(n-1,1) \to Imm^{sf}(n_0,n-n_0)$, see (\ref{102})
for the case $A_1 = \{\kappa_0=\kappa_M\}$.
Denote by $x_{d-1}$  the total obstruction of a compression of
the order $d-1$ for the element $x$. 
By Corollary \ref{totalobstr} the element $x_{d-1}$ coincides with the total obstruction of a
compression of the order  $(d-1)$ for the element $J^{sf}_{n_0}(x)$, $n_0 \ge d-1$.

Let us consider the total obstruction $x_{n_1} \in
Imm^{sf}(n_1,n-n_1)$ of the retraction of the order $n_1$ for the
element $x$, we have $x_{n_1}=J^{sf}_{n_1}(x)$. 
By Proposition $\ref{38}$, there exists an
element $y_{n_1} \in Imm^{sf;\{\kappa_1\}}(n_1,n-n_1)$, such that
the image of the element $y_{n_1}$ by the transfer homomorphism 
equals to $x_{n_1}$.

Let us consider the total obstruction $y_{n_2} \in
Imm^{sf;\kappa_1}(n_2,n-n_2)$ of a compression of the order $n_2$
for the element $y_{n_1}$. By the Proposition $\ref{38}$ there
exists an element  $z_{n_2} \in Imm^{sf; \kappa_1,
\kappa_2}(n_2,n-n_2)$ such that the element $r_2^!(z_{n_2})$ is the
total obstruction of a compression of the order  $n_2$ for the
element $y_{n_1}$. The element $r_{tot}(z_{n_2})=r_1 \circ
r_2(z_{n_2})=J^{sf}_{n_2}(x_{n_0})=x_{n_2}$ is the total
obstruction for a compression of the order  $n_2$ of the element
$x_{n_1}$. The same element $x_{n_2}$ is the total obstruction for
a compression of the order  $(n_2)$ for the initial element $x$.

By the induction we prove that the total obstruction of a
compression of the order  $(d-1)$ for the element $x$ is in the
image of the total transfer homomorphism of the multiplicity
 $\psi$, i.e. this total obstruction  equals to
 $r_{tot}(z)$,  $z \in Imm^{sf;A_{\psi}}(d-1,n-d+1)$. By Proposition
$\ref{totaltrans}$ we have  $r_{tot}(z)=0$. Therefore,
$x_{d-1}=0$. The Compression Theorem $\ref{comp}$ is proved. \qed
\[  \]

\section{Appendix: a self-contained proof of the Cohen's immersion Theorem \cite{C} in a special case}

In the proof of Proposition $\ref{38}$ we used the following fact.
Denote $b=b(k)=2^{2k}$.
Assume that $n=2^l-2$, $l \ge b(k)$. 

\begin{theorem}\label{39}
Let  a cobordism class $x \in
Imm^{sf,A_r}(b(k)-2,n-b(k)+2)$ is represented by a skew-framed immersion $(f: M^{b(k)-2} \looparrowright \R^n,\Xi_M)$,
$\dim(M^{b(k)-2})=b(k)-2$, where the manifold $M^{b(k)-2}$ is
equipped with a mapping $(\ref{lam})$.
Assume that the element $x$ belongs to the kernel of the forgetful homomorphism
\begin{eqnarray}\label{Omega}
Imm^{sf,A_r}(b(k)-2,n-b(k)+2) \to H_{b(k)-2}(\prod_{i=0}^r \RP^{\infty}_i).
\end{eqnarray}
Then in the regular cobordism class $[x]$ there exists an element for which the manifold $M^{b(k)}$ admits an immersion into the Euclidean
space $\R^{2b(k)-3-k}$ with a non-degenerate skew cross section given by a linear bundle $\kappa_M$.
\end{theorem}



\subsection*{Remark}
 Theorem $\ref{39}$
has an alternative proof, because it is a corollary of the R.Cohen's Immersion Theorem [C]. 
\[  \]
 
The main step of the proof of Theorem $\ref{39}$ is the following lemma.

\begin{lemma}\label{Cohen}
Let $(N^{j}, \Xi_N)$, $j=0, \dots, k-1$, be a framed (in particular, oriented) manifold, equipped with a mapping
$$\lambda_N: N^j \to \RP^{\infty}_0 \times \prod_{i=1}^r \RP^{\infty}(i).$$
Let us assume that the Hurewicz image
\begin{eqnarray}\label{Hur}
(\lambda_N)_{\ast}([N]) \in H_{j}(\prod_{i=1}^r \RP^{\infty}_i;\Z)
\end{eqnarray}
 with integer coefficients 
 is trivial. Then there exists a skew-boundary $(W^{k},\Psi_W, \lambda_W)$ in codimension $b(k)-1$,
 $\partial(W^{k}_W) = N^{k-1}$, $\Psi_W \vert_{\partial W^{k}} = \Xi_N$, $\lambda_W \vert {\partial W=N^k} = \lambda_N$, where
$$\lambda_W: W^{k} \to \RP^{\infty}_0 \times \prod_{i=1}^r \RP^{\infty}(i), $$
 and $\kappa_W = w_1(\Psi_W)$ coincides with the projection of the mapping $\lambda_W$ on the factor $\RP^{\infty}_0$.
\end{lemma}
\[  \]

\subsubsection*{A sketch of the proof of Lemma $\ref{Cohen}$}
 To prove Lemma  $\ref{Cohen}$ let us consider a spectral sequence of the Atiyah–Hirzebruch type  for the
cobordism group
of framed immersions of the dimensions $0, \dots, k-1$ (the stable homotopy group) of the space $\prod_{i=1}^r \RP^{\infty}_i$.
 Let us consider the Atiyah–Hirzebruch spectral sequence  for
cobordism groups  
of skew-framed immersions of dimensions $0, \dots, k-1$ in the codimension $b(k)-1$ of the space $\RP^{\infty}_{0} \times \prod_{i=1}^r \RP^{\infty}_i$.
There is a natural mapping of the first spectral sequence to the second spectral sequence.
By the main result of [A-E] all higher coefficients  in $E_2$-terms of the kernel  are trivial. Therefore the kernel of $E^{\infty}$--therm is totally described by
the Hurewicz image $(\ref{Hur})$ of a cobordism class of a corresponding mapping of a framed manifold. This Hurewicz image is trivial by the assumption. Therefore an arbitrary framed manifold $(N^{k}, \Xi_N)$ is a skew-framed boundary in codimension $b(k)-1$.
 Lemma Lemma $\ref{Cohen}$ is proved. \qed

\subsubsection*{A sketch of the proof of Theorem $\ref{39}$}

By the dimensional assumption, the normal bundle
 $\nu_M$ over $M$ is isomorphic to the Whitney sum  $Cb(k)\kappa_M$, where $\kappa_M$
 is a given line bundle over $M^{b(k)-2}$, where we may assume for a simplicity
$b(k)=2^{2k}$, $C$ is a positive integer. Let us calculate 
 the Koschorke construction  \cite{K} to prove that the normal bundle $\nu_M=Cb(k)\kappa_M$ over $M^{b(k)-2}$ 
admits a fibrewised monomorphism of the bundle  
\begin{eqnarray}\label{muM}
\mu_M=(C-1)b(k)\varepsilon \oplus  k\varepsilon \oplus \kappa_M .
\end{eqnarray}
 We will speak below about $\mu$-section of the
bundle $\nu_M$ for short.
The construction is given by the induction over the index $j=0, \dots, k$.

Let us denote that the bundle
$\nu_M$ admits a regular family of  $2k+(C-1)b(k)$ sections. 
Denote by 
$$I \subset \nu_M$$ 
the orthohonal complement to the subbundle $(C-1)b(k)\varepsilon  \subset \nu_M$, which consists of first regular
 $(C-1)b(k)$ $\mu$-sections, see (\ref{muM}).  
At the $j$-th induction step, $j=0, \dots, k$, let us assume that there exists
a morphism (which is not a fibrewised monomorphism)   $\rho_{j}$ of the subbundle  $(2k-j+1)\varepsilon \oplus \kappa_M \subset \mu_M $, which consists of the last terms in (\ref{muM}) into the subbundle $I \subset \nu_M$, such that the following two conditions are satisfied:

--1. Denote by  $\hat \rho_{j}$ the restriction of the morphism    $\rho_{j}$ on the subbundle 
$(j+1) \varepsilon \oplus \kappa_M \subset (2k-j+1)\varepsilon \oplus \kappa_M$, the restriction
should be a regular morphism, i.e. a fibrewised monomorphism. 

--2.  Denote by  $\tilde \rho_{j}$ the restriction of the morphism    $\rho_{j}$ on the subbundle 
$(2k-j+1)\varepsilon \subset (2k-j+1)\varepsilon \oplus \kappa_M$, the restriction
should be a regular morphism. 

The base step of the induction is provided by the condition, because the element $x$ in the kernel of  $(\ref{Omega})$.

Let us prove a step
$j \mapsto j+1$ of the induction. Consider the morphism $\rho_{j}$ and denote the restriction of this morphism to the subbundle 
$(j+2) \varepsilon \oplus \kappa_M \subset (2k-j+1)\varepsilon \oplus \kappa_M$ by $\hat \rho^{\circ}_j$ (one section more then in the family 
$\hat \rho_j$). The morphism $\hat \rho^{\circ}_j$, generally speaking, is not regular. 
Let us denote by 
  $N^{j+1} \subset M^{b(k)-2}$ the singular submanifold of this morphism.  Conditions imply
 that the restriction of the class  $\kappa_M$ on
$N^{j+1}$ is trivial and that $N_{j+1}$ is framed in  $M^{b(k)-2}$ (and, therefore, 
$N^{j+1}$ is a framed manifold). Let us prove that this framed manifiold satisfied the conditions of Lemma  $\ref{Cohen}$.

Let us prove that the  Hurewicz image $(\ref{Hur})$ is trivial.
There are the two cases: 

--a. $j$ is even; 

--b. $j$ is odd.

Let us consider the case a. Take the collection of section $\rho_{j+2}$, which is restricted to the subbundle
$(j+3) \varepsilon \oplus \kappa_M \subset (2k-j+1)\varepsilon \oplus \kappa_M$,
denote this restriction by $\hat \rho^{\circ \circ}_{j}$
(extra two sections with respect to $\hat \rho_j$). 
Singularities of the morphism $\hat \rho^{\circ \circ}_{j}$ is denoted by  $K^{j+2}$--this is a submanifold of non-regular sections. It is easy to calculate: 
$w_1(K^{j+2})=\kappa_M \vert_{K}$.  The condition that the total obstruction $x$ by the homomorphism $(\ref{Omega})$ is oriented, implies that  $K^{j+2}$, as an element in the group $H_{j+2}(\RP^{\infty}_0 \times\prod_{i=1}^r \RP^{\infty}_i;\Z^{tw})$, in which twisted coefficients is considered with respect to the first term in the Cartesian product and are given by the class $\kappa_M$, is a boundary. Therefore, the homology class  $[K^{j+2}] \cap \kappa_M = [L^{j+1}]$, which is considered as an element in the group  $H_{j+1}(\RP^{\infty}_0 \times\prod_{i=1}^r \RP^{\infty}_i;\Z)$, is a boundary. It is proved that in the case --a Lemma $\ref{Cohen}$ should be applied.

Let us consider the case --b. Consider the family
 $\rho_{j}^{\circ}$, which contains a regular sub-family 
$\rho_{j}$. Because a restriction  $\kappa_M$ on the singular manifold  $N^{j+1}$ is trivial, the
tensor product of the family  $\rho_{j}^{\circ}$ with the line bundle  $\kappa_M$ is degenerated. The family 
$\rho_{j}^{\circ} \otimes \kappa_M$ determines a morphism of the bundle  $(C-1)b(k)\kappa_M \oplus (j+1)\kappa_M \oplus \varepsilon$ into the bundle  $Cb(k)\varepsilon - b(k)\kappa_M$.
Because  $j+1$ is even, the image of the total obstruction $(\ref{Hur})$ is oriented, and 
we may use the condition that  $x$ in the kernel  of $(\ref{Omega})$. 
It is proved in the case --b that conditions of Lemma  $\ref{Cohen}$ are satisfied.

The obstruction of the existence of a regular morphism $\hat \rho_{j+1}$
with the trivial kernel over the singular manifold is given by a  cobordism class of a mapping $\lambda_N: N^j \to \prod_{i=1}^r \RP^{\infty}_i$ of a framed $j$-dimensional manifold, which is considered as a  skew-framed manifold in codimension $b(k)-1$ with local twisted coefficients system,
associated with $\kappa_M$  (note that $\kappa_M \vert_{N^{j+1}}$ is trivial, but in the regular cobordism class this property is not assumed). 

We proved that there exist $\rho_{j+1}$, which satisfied Condition --1. Let us prove that  $\rho_{j+1}$ can be taken with the Condition  2. By the construction,
a homotopy of the section $\hat \rho_{j}^{\circ}$ into a section  $\hat \rho_{j+1}$ has
is a  singular manifold  with a trivial kernel fibration.

Evidently, there exists a homotopy
 $\tilde \rho_j$ в $\tilde \rho'_j$ of the sub-family of sections $(2k-j+1)\varepsilon$, 
 which has singularities with the trivial kernel, moreover, the restriction  $\tilde \rho'_j$ on the subfamily $(j+1)\varepsilon$ coincids with the restriction  $\hat \rho_{j+1}$ on the family. Therefore, there exist a monomorphism
 $\rho'_{j+1}$ of the subbundle $(2k-j+1)\varepsilon \oplus \kappa_M$ (one section more then 
$\rho_{j+1}$, for which Condition  --1, and Condition --2  are satisfied and the following condition is satisfied: the restriction
$\rho'_{j+1}$ on the subbundle   $(2k-j+1)\varepsilon \subset (2k-j+1)\varepsilon \oplus \kappa_M$ has a trivial kernels fibration. By general position arguments, the restriction
 $\rho'_{j+1}$ on this trivial kernels fibration is regular.
 Let us restrict the morphism
$\rho'_{j+1}$ on the subbundle  $(2k-j)\varepsilon \subset (2k-j+1)\varepsilon \oplus \kappa_M$.
There existe a small generic deformation
 $\rho'_{j+1}\vert_{(2k-j)\varepsilon} \mapsto \tilde \rho_{j+1}$, into a regular morphism
 $\tilde \rho_{j+1}$. This proves Condition 2. The step of the induction is proved.

By the last step of the induction, there exists a regular section $\rho_{k-1}$.
Theorem $\ref{39}$ is proved.  \qed


Ахметьев Петр Михайлович  
 pmakhmet@mail.ru

\end{document}